\documentclass[10pt]{article}

\usepackage{amsmath,amsthm}

\newtheorem{thm}{Theorem}[section]
\newtheorem{cor}[thm]{Corollary}
\newtheorem{prop}{Proposition}
\newtheorem{lem}[thm]{Lemma}

\theoremstyle{remark}

\theoremstyle{plain}

\newtheoremstyle{note}
  {3pt}
  {3pt}
  {}
  {}
  {\itshape}
  {:}
  {.5em}
  {}

\theoremstyle{note}

\newtheoremstyle{citing}
  {3pt}
  {3pt}
  {\itshape}
  {}
  {\bfseries}
  {.}
  {.5em}
  {\thmnote{#3}}

\theoremstyle{citing}

\newtheoremstyle{break}
  {9pt}
  {9pt}
  {\itshape}
  {}
  {\bfseries}
  {.}
  {\newline}
  {}

\theoremstyle{break}

\theoremstyle{exercise}

\swapnumbers
\theoremstyle{plain}

\let\lvert=|\let\rvert=|

\begin{document}

\begin{center} 
\bf{Some Irrational Generalised Moonshine from Orbifolds}
\end{center}

\begin{center} 
Rossen Ivanov$^{1,2}$ and Michael Tuite$^{1,3}$\\ 

$^{1}$Department of Mathematical Physics,\\
National University of Ireland, Galway, Ireland\\

$^{2}$Institute for Nuclear Research and Nuclear Energy,\\
72 Tzarigradsko shosse, 1784 Sofia, Bulgaria\\

$^{3}$Dublin Institute for Advanced Studies,\\
10 Burlington Road, Dublin 4, Ireland\\

email: rossen.ivanov@nuigalway.ie, michael.tuite@nuigalway.ie
\end{center}

\begin{abstract}
We verify the Generalised Moonshine conjectures for some irrational modular
functions for the Monster centralisers related to the Harada-Norton, Held, $M_{12}$ and L$_{3}$(3) simple groups based on certain orbifolding
constraints. We find explicitly the fixing groups of the hauptmoduls arising in each case.

PACS: 11.25.Hf, 02.10.De, 02.20.Bb

Key Words: Conformal Fields, Modular Groups, Moonshine, Orbifolds.
\end{abstract}

\section{\protect\smallskip Introduction}

The Moonshine Module \cite{FLM1},\cite{FLM2} whose automorphism group is the
Monster finite sporadic group $\bf{M}$, is an orbifold Meromorphic
Conformal Field Theory (MCFT) constructed by orbifolding the Leech lattice
MCFT with respect to the group generated by a reflection involution. The
Moonshine Module partition function is the classical elliptic $J$ function
which is a hauptmodul for the genus zero modular group $SL(2,\bf{Z})$. The
Moonshine Module is believed to be the unique MCFT with this partition
function \cite{FLM2}. Orbifolding the Moonshine Module with respect to the
group generated by an element $g\in $ $\bf{M}$ leads naturally to the
notion of the orbifold partition function known as the Thompson series $%
T_{g} $. Monstrous Moonshine is mainly concerned with the property,
conjectured by Conway and Norton \cite{CN} and subsequently proved by
Borcherds \cite{B}, that each Thompson series $T_{g}$ is a hauptmodul for
some genus zero fixing modular group. Assuming the uniqueness of the
Moonshine Module, this genus zero property is also believed to be equivalent
to the following statement \cite{T2}: the only orbifold MCFT that can arise
by orbifolding with respect to $g$ is either the Moonshine Module itself
(for $g$ belonging to a so-called Fricke Monster conjugacy class) or the
Leech lattice MCFT (for $g$ belonging to a non-Fricke Monster conjugacy
class).

The Generalised Moonshine conjecture of Norton \cite{N1} is concerned with
modular functions associated with a commuting pair $g,h\in $ $\bf{M}$ and
asserts that each such modular function known as a Generalised Moonshine
Function (GMF) is either constant or is a hauptmodul for some genus zero
fixing group. No\ extension of the Borcherds' approach to Monstrous
Moonshine has yet been shown to be possible for Generalised Moonshine. The
most natural setting for these conjectures is to consider orbifoldings of
the Moonshine Module with respect to the abelian group $\langle g,h\rangle $
generated by $g,h$ \cite{T3},\cite{IT}. In \cite{IT} we considered the case
where $g$ is of prime order $p=2,3,5$ and $7$ and is of Fricke type and $h$
is of order $pk$ for $k=1$ or $k$ prime. There we confirm Norton's
conjecture for modular functions with rational coefficients in these cases
by considering orbifold modular properties and some consistency conditions
arising from the orbifolding procedure leading to constraints on the
possible Monster conjugacy classes to which the elements of $\langle
g,h\rangle \,$ may belong. In the present paper we extend the approach based
on the restrictions coming from the orbifolding of the Moonshine Module and
demonstrate the hauptmodul property for GMFs with irrational coefficients
where $g$ and $h$ are of the same prime order.

We begin in Section 2 with a brief review of some properties of Meromorphic
Conformal Field Theory (MCFT), the Moonshine Module and Generalised
Moonshine Functions. In Section 3 we discuss the general properties for GMFs
with irrational coefficients when $g$ is of prime order. We then state two
theorems about the conjugation properties of the centraliser elements, and
one theorem concerning constraints that arise from the consistency of
orbifolding the Moonshine Module with respect to $\langle g,h\rangle $ under
specific choices of generators. We then analyse the GMFs with irrational
coefficients in the cases where $g$ is Fricke of prime order $p=5,7,11$ and $%
13$ and $h$ is of same order $p$. This analysis in part relies on properties
of the characters of the centralisers of the Monster related to the
Harada-Norton, Held, Mathieu and $\mathrm{L}_{3}\mathrm{(3)}$ simple groups.
We give a comprehensive analysis of the possible singularity structure of
GMFs for the cases under consideration. In each case we demonstrate that all
singularities of the GMF can be identified under some genus zero fixing
group for which the GMF is a hauptmodul. Thus we verify the Generalised
Moonshine conjecture in these cases.

\section{Generalised Moonshine}

\subsection{\textbf{\protect\smallskip }Self-Dual Meromorphic Conformal
Field Theories}

\textbf{\ }Let $\mathcal{H}$ denote the Hilbert space of a Self-Dual
Meromorphic Conformal Field Theory (MCFT) \cite{Go} of central charge $24$.
The characteristic function (or genus one partition function) for $\mathcal{H%
}$ is given by $Z(\tau )=Tr_{\mathcal{H}}(q^{L_{0}-1}),$ $q=e^{2\pi i\tau }$%
, where $\tau \in \bf{H}$, the upper half complex plane, is the usual
elliptic modular parameter. $\mathcal{H}$ has integral grading (the
conformal weight) with respect to $L_{0}$, the Virasoro level operator so
that $Z(\tau )$ $\,$is invariant under $T:\tau \rightarrow \tau +1$.
Furthermore, self-duality of the MCFT implies invariance under $S:\tau
\rightarrow -1/\tau $ so that $Z(\tau )$ is invariant under the modular
group ${SL}(2,\bf{Z})$ generated by $S,T$. Hence $Z(\tau )$ is given up to
an additive constant by $J(\tau )$, the hauptmodul for ${SL}(2,\bf{Z})$ 
\cite{Se} i.e. ${Z(\tau )}={J(\tau )+N_{0}}$, where $J$ has $q$ expansion 
\begin{equation}
J(\tau )=\frac{E_{4}^{3}(\tau )}{{\eta ^{24}(\tau )}}-744=\frac{1}{q}%
+0+196884q+21493760q^{2}+...  \label{J}
\end{equation}
and where $\eta (\tau )=q^{1/24}\prod_{n>0}(1-q^{n})$ is the Dedekind eta
function, $E_{n}(\tau )$ is the Eisenstein modular form of weight $n$ \cite
{Se}, $N_{0}$ is the number of conformal weight 1 operators in $\mathcal{H}$%
. For the Leech lattice MCFT $N_{0}=24$ and for the FLM Moonshine Module $%
\mathcal{V}^{\natural }$ with Hilbert space $\mathcal{H}^{\natural }$ \cite
{FLM1,FLM2} $N_{0}=0$. This means that the latter does not contain conformal
dimension 1 operators i.e. there is an absence of the usual Kac-Moody
symmetry.

\subsection{The Moonshine Module and Monstrous Moonshine}

The Monster group $\bf{M}$, the largest finite sporadic simple group, is
the automorphism group of the Moonshine Module \cite{FLM1,FLM2}. The
Thompson series $T_{g}(\tau )$ for each $g\in \bf{M}$ is defined by 
\begin{equation}
T_{g}(\tau )\equiv \mathrm{Tr}_{\mathcal{H}^{\natural }}(gq^{L_{0}-1})
\label{Thompson_g}
\end{equation}
with $q$ expansion with coefficients determined by (reducible) characters of 
$\bf{M}$. The Thompson series for the identity element is $J(\tau )\,$of (%
\ref{J}), which is the hauptmodul for the genus zero modular group $\mathrm{%
SL}(2,\bf{Z})$ as already stated.

Conway and Norton \cite{CN} conjectured and Borcherds \cite{B} proved that $%
T_{g}(\tau )$\ is the hauptmodul for some genus zero fixing modular group $%
\Gamma _{g}$. This remarkable property is known as Monstrous Moonshine. In
general, for $o(g)=n$, $T_{g}(\tau )$ is found to be $\Gamma _{0}(n)$
invariant up to $m^{\mathrm{th}}$ roots of unity where 
\[
\Gamma _{0}(n)\equiv \left\{ \left( 
\begin{array}{cc}
a & b \\ 
c & d
\end{array}
\right) \in \Gamma ,\ c=0\ {\rm mod}\ n\right\} , 
\]
and where $m$ is an integer with $m|n$ and $m|24$. $g$ is said to be a
normal element of $\bf{M}$ if and only if $m=1$, otherwise $g$ is said to
be anomalous. $T_{g}(\tau )$ is fixed by some $\Gamma _{g}\supseteq \Gamma
_{0}(N)$ which is contained in the normalizer of $\Gamma _{0}(N)$ in $SL(2,%
\bf{R})$ where $N=nm$ \cite{CN}. This normalizer contains the Fricke
involution $W_{N}:\tau \rightarrow -1/N\tau $. All classes of $\bf{M}$ can
therefore be divided into Fricke and non-Fricke type according to whether or
not $T_{g}(\tau )$ is invariant under the Fricke involution. There are a
total of 51 non-Fricke classes of which 38 are normal and there are a total
of 120 Fricke classes of which 82 are normal. The genus zero properties of
Monstrous Moonshine can also be understood using constraints arising from
the possible orbifolding of $\mathcal{V}^{\natural }$ with respect to the
group generated by $g$ \cite{T1}, \cite{T2}.

\subsection{Properties of Generalised Moonshine Functions}

We now consider Generalised Moonshine Functions (GMFs) which are generalised
Thompson series depending on two commuting Monster elements of the form 
\begin{equation}
Z\left[ 
\begin{array}{c}
h \\ 
g
\end{array}
\right] (\tau )\equiv \mathrm{Tr}_{\mathcal{H}_{g}^{\natural
}}(hq^{L_{0}-1}),  \label{Z}
\end{equation}
for $h\in C_{g}$, the centralizer of $g$ in $\bf{M}$. $\mathcal{H}%
_{g}^{\natural }$ denotes the Hilbert space for a so called $g$-twisted
sector of $\mathcal{V}^{\natural }$ \cite{DLM2},\cite{IT}. The original
Thompson series (\ref{Thompson_g}) then corresponds to the untwisted sector,
i.e. 
\[
T_{g}(\tau )=Z\left[ 
\begin{array}{c}
g \\ 
1
\end{array}
\right] (\tau ).
\]
Each $h\in C_{g}$ has a central extension, acting on $\mathcal{H}%
_{g}^{\natural }$, which we also denote by $h$ \cite{IT}. Norton has
conjectured that \cite{N1}:

\noindent{\bf Generalised Moonshine.}{\it 
The GMF (\ref{Z}) is either constant or is a hauptmodul for some genus zero
fixing group $\Gamma _{h,g}$.
}

If $\langle g,h\rangle =\langle u\rangle $ for some $u\in \bf{M}$ then (\ref
{Z}) can always be transformed to a regular Thompson series (\ref{Thompson_g}%
) \cite{T1},\cite{DLM1},\cite{IT}. In particular this is always possible
when $o(g)\,$and $o(h)$ are coprime. In these cases, the genus zero property
for the GMF therefore follows from that for a regular Thompson series (\ref
{Thompson_g}). The GMFs with non-trivial genus zero behaviour then occur for 
$h\in C_{g}$ where $o(h)$ and $o(g)$ are not coprime.

In Section 3 we analyse GMFs with irrational coefficients in their $q$%
-expansion for $o(g)=o(h)=p$, prime, where $\langle g,h\rangle \simeq {\bf Z}%
_{p^{2}}$ (then $g,h$ are independent, $\,$i.e. $g^{A}\neq h^{B}$ for all $%
A,B=1,...,p-1$). We denote $\,$the fixing group for $Z\left[ 
\begin{array}{c}
h \\ 
g
\end{array}
\right] (p\tau )\,$ by $\tilde{\Gamma}_{h,g}$ which is obviously conjugate
to $\Gamma _{h,g}$ where 
\[
\tilde{\Gamma}_{h,g}=\left\{ \widetilde{\gamma }|\widetilde{\gamma }=\theta
_{p}\gamma \theta _{p}^{-1},\;\gamma \in \Gamma _{h,g},\;\theta _{p}\equiv
\left( 
\begin{array}{cc}
1 & 0 \\ 
0 & p
\end{array}
\right) \right\} . 
\]
Various properties of (\ref{Z}) arising from orbifold considerations can be
summarised as follows \cite{IT}:

(i) When all elements of $\langle g,h\rangle $ $\,$are normal elements of $%
\bf{M}$ then for $\gamma =\left( 
\begin{array}{cc}
a & b \\ 
c & d
\end{array}
\right) \in \Gamma \equiv SL(2,{\bf Z})$

\begin{equation}
Z\left[ 
\begin{array}{c}
h \\ 
g
\end{array}
\right] (\gamma \tau )=Z\left[ 
\begin{array}{c}
h^{d}g^{-b} \\ 
h^{-c}g^{a}
\end{array}
\right] (\tau ),\;\gamma \tau =\frac{a\tau +b}{c\tau +d}  \label{SLT}
\end{equation}
Hence $\Gamma (p)\subseteq \Gamma _{h,g}$, (in fact $\Gamma (p)\triangleleft
\Gamma _{h,g}$) with 
\[
\Gamma (p)\equiv \left\{ \left( 
\begin{array}{cc}
a & b \\ 
c & d
\end{array}
\right) \in \Gamma ,\ \left( 
\begin{array}{cc}
a & b \\ 
c & d
\end{array}
\right) =\left( 
\begin{array}{cc}
1\ {\rm mod}\ p & 0\ {\rm mod}\ p\  \\ 
0\ {\rm mod}\ p & 1\ {\rm mod}\ p
\end{array}
\right) \right\} 
\]
In particular, since $\gamma $ and $-\gamma $ act equally we have

\begin{equation}
Z\left[ 
\begin{array}{c}
h \\ 
g
\end{array}
\right] =Z\left[ 
\begin{array}{c}
h^{-1} \\ 
g^{-1}
\end{array}
\right] ,  \label{Z1}
\end{equation}
which property is known as charge conjugation invariance.

(ii) For a normal non-Fricke element $g$ of order $p$ 
\begin{equation}
Z\left[ 
\begin{array}{c}
h \\ 
g
\end{array}
\right] (\tau )=\text{const}+O(q^{1/p}).
\end{equation}
i.e. $\mathcal{H}_{g}^{\natural }$ has zero vacuum energy and the GMF is
nonsingular at $q=0$ \cite{T1},\cite{T3}$.$

The value of (\ref{Z}) at any parabolic cusp $a/c$ with $(a,c)=1$ is
determined by the vacuum energy of the $g^{a}h^{-c}$ twisted sector from (%
\ref{SLT}). If all elements of $\langle g,h\rangle $ are non-Fricke then
there are no singular cusps so that $Z\left[ 
\begin{array}{c}
h \\ 
g
\end{array}
\right] $ is holomorphic on ${\bf H}/\Gamma (p)$ and hence is constant. This
accounts for the constant GMFs referred to in the Generalised Moonshine
Conjecture above. We therefore assume from now on that at least one element
of $\langle g,h\rangle $ $\,$is Fricke which we chose to be $g$ without loss
of generality.

(iii) For a normal Fricke element $g$ of order $p$ the vacuum energy of $%
\mathcal{H}_{g}^{\natural }$ is $-1/p$ and the GMF is singular at $q=0$ \cite
{IT}: 
\begin{equation}
Z\left[ 
\begin{array}{c}
h \\ 
g
\end{array}
\right] (\tau )=\phi _{g}(h)[q^{-1/p}+0+{\sum_{s=1}^{\infty }}%
a_{g,s}(h)q^{s/p}].  \label{Z(h,Fricke)}
\end{equation}
Here $\phi _{g}(h)$ is a $p$-th root of unity describing the central
extension of the action of $h\in C_{g}$ on $\mathcal{H}_{g}^{\natural }$ 
\cite{IT}. It is further assumed that $\phi _{g}$ can be chosen so that \cite
{IT} 
\begin{equation}
\phi _{g}(g^{a}h^{b})=\phi _{g}(g)^{a}\phi _{g}(h)^{b},\;\phi _{g}(1)=1.
\label{phi g a h b}
\end{equation}
Then by considering a $T$ transformation in (\ref{SLT}) with $h=1$ we have

\begin{equation}
\phi _{g}(g)=\omega _{p}\equiv \exp (2\pi i/p)  \label{FA}
\end{equation}
Let $G_{g}\equiv C_{g}/\langle g\rangle $. If $C_{g}=\langle g\rangle \times
G_{g}$ then $h$ and $gh$ are not elements of the same $C_{g}$ conjugacy
class and we assume that $\phi _{g}\,$ can be chosen so that \cite{IT} 
\begin{equation}
\phi _{g}(xhx^{-1})=\phi _{g}(h)  \label{phi_class_fun}
\end{equation}
for all $x\in C_{g}$ . By inspection from the ATLAS \cite{CCNPW} $h\in G_{g}$
is conjugate to $h^{a}$, for some $a\neq 0\ {\rm mod}\ p$, so that (\ref
{phi g a h b}) implies 
\begin{equation}
\phi _{g}(h)=1.  \label{phi_g(h)}
\end{equation}
The coefficient $a_{f,s}(h)$ is called a head character for the given GMF
and is a reducible character of $G_{g}$ \cite{DLM1},\cite{IT}. If $%
g^{a}h^{-c}\,$is Fricke then due to (\ref{Z(h,Fricke)}) and (\ref{SLT}) the
GMF is singular at $a/c$ with residue $\phi _{g^{a}h^{-c}}(h^{d}g^{-b})$.
The GMF (\ref{Z}) is holomorphic at all other points on $\bf{H}$. Once
these singularities are known, then (\ref{Z}) can be analysed to check
whether it is constant or is a hauptmodul for an appropriate genus zero
modular group. The singularities and residues of $Z\left[ 
\begin{array}{c}
h \\ 
g
\end{array}
\right] $ are restricted by certain orbifolding constraints discussed in 
\cite{IT}.

(iv) Given the uniqueness of the twisted sectors for $\mathcal{H}^{\natural
} $, under conjugation by any element $x\in \bf{M}$ then $x(\mathcal{H}%
_{g}^{\natural })x^{-1}$ is isomorphic to $\mathcal{H}_{xgx^{-1}}^{\natural
} $ so that if $C_{g}=\langle g\rangle \times G_{g}$ with $h\in G_{g}$ \cite
{IT}

\begin{equation}
Z\left[ 
\begin{array}{c}
h \\ 
g
\end{array}
\right] =Z\left[ 
\begin{array}{c}
xhx^{-1} \\ 
xgx^{-1}
\end{array}
\right] .  \label{theta}
\end{equation}

\section{Irrational GMFs for $g=p+$}

\subsection{Monster centralisers and irrational characters}

We now briefly describe all irrational GMFs (\ref{Z}) where $g$ is Fricke of
prime order $p$ i.e. $g=p+$ in Conway-Norton notation \cite{CN} and $h$ of
order $p$. A GMF is said to be irrational if it has at least one irrational
head character (\ref{Z(h,Fricke)}). A character is said to be irrational if
it is irrational for at least one conjugacy class. Otherwise, it is said to
be rational. Rational cases are discussed in \cite{IT}.

From the ATLAS \cite{CCNPW} we see that irrational characters occur for $p=2$%
, $3$, $5$, $7$, $11$ and $13\,$. From the explicit calculations for the
head characters of $C_{2+}=2.\mathrm{B}$ ($p=2$) and $C_{3+}=3.\mathrm{Fi}$ (%
$p=3$) \cite{IT} we observe that there are no irrational head characters for 
$p=2$ and that there is only one pair of conjugate irrational head
characters for $p=3$, with $h$ of order $54$. For $p=5$, $7$, $11$ and $13$ $%
\,$the centraliser $C_{p+}=\langle p+\rangle \times G_{p+}$ for simple group 
$G_{p+}$ is as follows:

\begin{tabular}{|l|l|l|l|}
\hline
$g=p+$ & $C_{p+}$ & $G_{p+}$ & Name \\ \hline
$5+$ & $\mathrm{5\times HN}$ & $\mathrm{HN}$ & Harada-Norton \\ \hline
$7+$ & $\mathrm{7\times He}$ & $\mathrm{He}$ & Held \\ \hline
$11+$ & $\mathrm{11\times M}_{12}$ & $\mathrm{M}_{12}$ & Mathieu \\ \hline
$13+$ & $\mathrm{13\times L}_{3}\mathrm{(3)}$ & $\mathrm{L}_{3}\mathrm{(3)}$
&  \\ \hline
\end{tabular}

Table 1. $p+$ centralisers for $p=5,7,11,13$

We will focus now our analysis on the groups from Table 1. For $h\in G_{g}$
from (\ref{phi_g(h)}) the GMF (\ref{Z(h,Fricke)}) has leading $q$ expansion in so-called
Normalised Function Form \cite{FMN}: 
\begin{equation}
Z\left[ 
\begin{array}{c}
h \\ 
g
\end{array}
\right] (\tau )=\frac{1}{q^{1/p}}+0+O(q^{1/p}).
\end{equation}
This is sufficient to ensure that $g,h$ are independent generators of $%
\langle g,h\rangle $ since otherwise $h^{B}=g^{A}$ for some $A\neq 0\ {\rm mod}  p$ implies the contradictory relation $(\phi _{g}(h))^{B}=(\phi
_{g}(g))^{A}=\omega _{p}^{A}$ from (\ref{FA}). We make the following useful
observations concerning the irreducible characters for the groups from Table
1, which can be checked by inspecting the appropriate ATLAS character tables 
\cite{CCNPW}. In all cases irrational characters can only occur when $o(h)$
is divisible by $p$. When $p=5,7,11$ any irrationality for an irreducible
character $\chi $ is quadratic i.e. for $h\in G_{p+}$, $\chi (h)=a\pm \sqrt{b}$%
, $\overline{\chi }(h)\equiv \chi (h^{d})=a\mp\sqrt{b}$ ($d=2$ when $p=5$ and $%
d=-1$ when $p=7,11$, see below), for some $a,b\in \bf{Q}$ where $b=0$ for
rational $\chi (h)$. The \textrm{HN }group contains one pair of conjugate
irrational classes of orders $5$, $15$, $20$, $25$ and $30$, and two pairs
of conjugate irrational classes of order $10$. The \textrm{He }group
contains one pair of conjugate irrational classes of orders $21$ and $28$,
and two pairs of conjugate irrational classes of orders $7$ and $14$. The $%
\mathrm{M}_{12}$\textrm{\ }group contains only one pair of conjugate
irrational classes of order $11$. The only non-quadratic irrationalities
occur when $p=13$ for $\mathrm{L}_{3}\mathrm{(3)}$ which has four conjugate
irrational classes of order $13$. 

\begin{thm}
Consider $g=p+$ and $h\in G_{p+}$ with $o(h)$ divisible by $p$ for $%
p=5,7,11,13$. Then for integers $a$, $b\neq 0$ ${\rm mod}\ p$, $Z\left[ 
\begin{array}{c}
h^{a} \\ 
g^{b}
\end{array}
\right] $ is either equal to $Z\left[ 
\begin{array}{c}
h \\ 
g
\end{array}
\right] $ or one of its algebraic conjugates.
\end{thm}

\smallskip \textbf{Proof.} For $g=p+$ we have $xgx^{-1}=g^{b}$ for some $%
x\in \bf{M}$. For any irreducible character $\chi $ of $G_{g}$ and for all $%
h\in G_{g}$, $\chi ^{x}(h)\equiv \chi (xhx^{-1})$ $\,$is also an irreducible
character of the same dimension. Furthermore $xG_{g}x^{-1}=G_{g^{b}}=G_{g}$
implies that the number of elements in the conjugacy classes in $G_{g}$ for $%
h$ and $xhx^{-1}$ are equal. Since $h$ $\in G_{g}$, $h$ is conjugate to $%
h^{a}$ for some $a\neq 1\ {\rm mod}\ o(h)$ from an inspection of the
ATLAS character tables \cite{CCNPW}. Hence $xhx^{-1}$ is conjugate to $%
(xhx^{-1})^{a}$ and is also member of a $G_{g}$ class. Therefore $\phi
_{g^{b}}(xhx^{-1})=1$ also. For $p=5,7,11$ one can show that if $\chi
^{x}\neq \chi $, then $\chi ^{x}=\overline{\chi }$, the algebraic conjugate 
\cite{IT}. Thus $xhx^{-1}$ is conjugate to a power of $h$ in $G_{g}$. When $%
p=13$ there are four irrational characters that are algebraically conjugate.
Then $\chi ^{x}=\chi *2^{r}$, $r=0,1,2$ or $3$, where by definition $(\chi
*k)(h)\equiv \chi (h^{k})$ i.e. the operation $*k$ replaces every root of
unity in the character value, prime to $k$ by its $k^{\text{th}}$ power (see
below). This implies that $xhx^{-1}$ is either an element of the same class
of $G_{g}$ as $h$, or of a class algebraically conjugate to $h$ i.e. there
exists $A$ such that $xhx^{-1}\stackrel{G_{g}}{\sim }h^{A}$. Therefore from (%
\ref{theta}) in all cases there exists $A$ such that 
\begin{equation}
Z\left[ 
\begin{array}{c}
h \\ 
g
\end{array}
\right] =Z\left[ 
\begin{array}{c}
h^{A} \\ 
g^{b}
\end{array}
\right] .  \label{Delta like inv}
\end{equation}
The latter is either equal to $Z\left[ 
\begin{array}{c}
h^{a} \\ 
g^{b}
\end{array}
\right] $ or to one of its algebraic conjugates and hence the result
follows. $\qed $

From now on we will assume that $o(h)=p$, prime. Then the set of elements 
\begin{equation}
\mathcal{S}\equiv \left\{ h^{a}\in G_{g}\text{, }a=1,2,...,p-1\right\} ,
\label{h-a sets}
\end{equation}
with irrational characters can be divided into classes characterized by the
so-called Dirichlet character of $a$ as follows. Let $C_{p-1}\equiv \langle
*m\rangle $ be the cyclic group of order $p-1$ acting on $\mathcal{S}$
generated by 
\[
\ast m:h\rightarrow h^{m}\text{ }
\]
i.e. the Galois group for the polynomial $(x^{p}-1)/(x-1)$. All irrational
characters $\chi $ of $G_{g}$ are fixed by a cyclic subgroup $\langle
*(mN)\rangle $ of $C_{p-1}$, for some $N$. Then $\mathcal{S}$ can be divided
into classes with representatives 
\begin{equation}
\left\{ h^{a}\in \mathcal{S}|a\in C_{p-1}/\langle
*(mN)\rangle \equiv C_{N}\right\} ,  \label{representatives}
\end{equation}
characterised uniquely by a Dirichlet character of order $N$ defined as
follows. A Dirichlet character $\chi _{D}$ of the Galois group $C_{p-1}$ is
a homomorphism:

\begin{equation}
\chi _{D}:C_{p-1}\rightarrow {\bf C}^{*}.  \label{Dirichlet homomorphism }
\end{equation}
Then these characters form a group under pointwise multiplication which is
isomorphic to $C_{p-1}$, but not under any canonical map. The character of
order $N$, $\chi _{D}^{(N)}$ then gives an isomorphism 
\[
\chi _{D}^{(N)}:C_{p-1}/\langle *(mN)\rangle \rightarrow \langle \omega
_{N}\rangle 
\]
and thus characterises uniquely the representatives (\ref{representatives}).

For $p$ odd, the character group (\ref{Dirichlet homomorphism }) is cyclic
of even order $p-1$, hence it has an element of order $2$, $\chi _{D}^{(2)}$. 
This is the Legendre symbol: $\chi _{D}^{(2)}(a)\equiv \left( \frac{a}{p}%
\right) \equiv a^{(p-1)/2}$ mod $p=\pm 1$ which characterises (\ref
{representatives}) in the case of quadratically irrational characters for
which $N=2$. This applies to the characters of $G_{p+}$ for $p=5,7,11$.

There is a character of order $N=4$ precisely when $4|$ $(p-1)$ (e.g. $p=13$%
). There are then exactly two characters of order $4$, because a cyclic
group of order divisible by $4$ has exactly two elements of order $4$. These
two characters are inverse to each other and hence are complex conjugates of
each other. Their quotient will be the Legendre symbol. For $p=13$, we may
take $m=2$ as a generator of the Galois group $C_{12}$. The fixing group for
the irrational characters of $G_{13+}=\mathrm{L}_{3}\mathrm{(3)}$ is then $%
\langle *3\rangle $ \cite{CCNPW}. We get a Dirichlet character $\chi
_{D}^{(4)}$ of order $4$ by mapping $2$ to $\omega _{4}=i=\exp (\pi i/2)$.
The other element of order $4$ obtains by mapping $2$ to $-i$. If $a$ is in $%
C_{12}$ with $a=2^{r}$ mod 13 then $\chi _{D}^{(4)}(a)=\exp (r\pi i/2)$. Thus 
\[
\chi _{D}^{(4)}:\text{%
\begin{tabular}{l}
$1,3,9\rightarrow 1$ \\ 
$12,10,4\rightarrow -1$ \\ 
$5,2,6\rightarrow i$ \\ 
$8,11,7\rightarrow -i$%
\end{tabular}
} 
\]
So $1,12,5,8$ are representatives of $C_{12}/\langle *3\rangle $ which map
to the four different $4^{\text{th}}$ roots of unity. Thus the value of $%
\chi _{D}^{(4)}$ characterises each of the four conjugate irrational classes
of order $13$ in $\mathrm{L}_{3}\mathrm{(3)}$.

\subsection{The Genus Zero Property for Irrational GMFs for $g=p+$ and $%
o(h)=p$}

We now come to the main purpose of this paper which is to demonstrate the
genus zero property for Generalised Moonshine Functions (GMFs) (\ref{Z}) in
the irrational cases where $g=p+$ and $o(h)=p$. Such cases occur when $%
p=5,7,11,13.$ Our aim is to show that the fixing group $\Gamma _{h,g}$
permutes all of the singular points of the GMF. In each case we find that $%
\Gamma _{h,g}$ is a subgroup of $\Gamma =SL(2,{\bf Z})$ so that $\Gamma
_{h,g}/\Gamma (p)$ is a subgroup of $L_{2}(p)=\Gamma /\Gamma (p)$, the group
permuting all ${\bf H}/\Gamma (p)$ inequivalent cusps. Then all possible
singularities of the GMF at the cusps $\tau =a/c$ with $(a,c)=1$ (related to
the $g^{a}h^{-c}$ Fricke-twisted sector) are identified under $\Gamma _{h,g}$
and the corresponding GMF is a hauptmodul for a genus zero group.

In our analysis we use the following two theorems. The first one gives some
constraints following from the orbifolding of the Moonshine Module with
respect to the abelian group $\langle g,h\rangle $. See \cite{IT} for
details:

\begin{thm}
\label{theorem_proj1} Let $g,h\in \bf{M}$ be independent commuting elements
where both $g$ and $h$ are Fricke such that $\phi _{g}(h)=1$ and where all
elements of $\langle g,h\rangle \,$ are normal. Let $u$, $v$ be any
independent generators for $\langle g,h\rangle $. If $u$ is Fricke then
there is a unique $A\ {\rm mod} \ o(u)$ such that $u^{A}v$ is Fricke with $%
\phi _{u^{A}v}(u)=1$ and $o(u^{A}v)=o(v)$.
\end{thm}

The next theorem gives the relation between the fixing groups of hauptmoduls
with algebraically conjugate coefficients \cite{N2},\cite{FMN}:

\begin{thm}
Let $G$ be a fixing group of a hauptmodul such that $\Gamma (p)\subseteq G$.
Let $k$ be any integer, coprime to $p$. Choose coset representatives of $%
\Gamma (p)$ in $G$, $\left( 
\begin{array}{cc}
a_{i} & b_{i} \\ 
c_{i} & d_{i}
\end{array}
\right) $, such that $k|c_{i}$. Then the operation $*k$, which replaces
every $p^{\text{th}}$ root of unity by its $k^{\text{th}}$ power yields a
Hauptmodul with a fixing group 
\[
G*k=\left\langle \Gamma (p),\left\{ \left( 
\begin{array}{cc}
a_{i} & kb_{i} \\ 
c_{i}/k & d_{i}
\end{array}
\right) \right\} \right\rangle 
\]
which is independent of choices made.
\end{thm}

In terms of Generalised Moonshine the action of $*k$ on the GMF $Z\left[ 
\begin{array}{c}
h \\ 
g
\end{array}
\right] $ is to map it to $Z\left[ 
\begin{array}{c}
h^{k} \\ 
g
\end{array}
\right] $ because $*k$ takes a character of $h$ to that of $h^{k}$. Hence 
\begin{equation}
\Gamma _{h^{k},g}=(\Gamma _{h,g})*k  \label{Gamma*k f.gr.}
\end{equation}

We begin with a preliminary Lemma.

\begin{lem}
\label{delta-p lemma} For $p=5,7,11,13$ we have $\Gamma _{0}^{0}(p)\equiv
\langle \Gamma (p),\delta _{p}\rangle $ for $\delta _{p}$ of order $(p-1)/2$
in $\Gamma (p)$, where 
\[
\Gamma _{0}^{0}(p)\equiv \left\{ \left( 
\begin{array}{cc}
a & b \\ 
c & d
\end{array}
\right) \in \Gamma ,\ b=c=0\ {\rm mod}\ p\right\} 
\]
and $\delta _{5}=\left( 
\begin{array}{cc}
2 & 5 \\ 
5 & 13
\end{array}
\right) $, $\delta _{7}=\left( 
\begin{array}{cc}
2 & 7 \\ 
7 & 25
\end{array}
\right) $, $\delta _{11}=\left( 
\begin{array}{cc}
-40 & -11 \\ 
11 & 3
\end{array}
\right) $, $\delta _{13}=\left( 
\begin{array}{cc}
85 & 13 \\ 
13 & 2
\end{array}
\right) $.
\end{lem}

\textbf{Proof}: By inspection of the possible solutions for $a,d$ 
mod $p$ of $ad-bc=1$ where $b$, $c=0$ mod $p$ one can
check the validity of the statement. $\qed $

Note that $\Gamma _{0}(p^{2})$ and $\Gamma _{0}^{0}(p)$ are conjugate where $%
\Gamma _{0}(p^{2})=$ $\theta _{p}\Gamma _{0}^{0}(p)\theta _{p}^{-1}$.
Therefore $\Gamma _{0}(p^{2})$ invariance follows from invariance under $%
\widetilde{\Gamma }(p)$ and $\widetilde{\delta }_{p}$.

\subsubsection{Irrational GMFs for $g=5+$ and $o(h)=5$}

The \textrm{HN} group has one pair of quadratically irrational classes of
order $5$. Then $m=2$, $N=2$ in (\ref{representatives}). The two subsets of $%
\mathcal{S}$ (\ref{h-a sets}) conjugate in $G_{5+}\equiv \mathrm{HN}$ are
characterised by the Legendre symbol $\left( \frac{a}{5}\right) $: $h%
\stackrel{G_{5+}}{\sim }h^{4}$ (i.e. $h^{a}$ where $\left( \frac{a}{5}%
\right) =1$) and $h^{2}\stackrel{G_{5+}}{\sim }h^{3}$ ($h^{a}$ where $\left( 
\frac{a}{5}\right) =-1$) \cite{CCNPW}. Therefore $gh\stackrel{{\bf M}}{\sim }%
gh^{4}$ and $gh^{2}\stackrel{{ \bf M}}{\sim }gh^{3}$. Since $g_{0}\stackrel{%
{\bf M}}{\sim }g_{0}^{s}$ for $s\neq 0$ mod $5$ for any Monster
element $g_{0}$ of order $5$ we have the following two disjoint sets $%
\mathcal{S}_{1},\mathcal{S}_{2}$ of conjugate elements in ${ \bf M}$ defined
by

\begin{eqnarray*}
\mathcal{S}_{1} &\equiv &\{(gh^{n})^{s}|\left( \frac{n}{5}\right) =1,s\neq 0\ {\rm mod}\ 5\} \\
\mathcal{S}_{2} &\equiv &\{(gh^{n})^{s}|\left( \frac{n}{5}\right) =-1,s\neq 0\ {\rm mod}\ 5\}
\end{eqnarray*}
i.e. $\mathcal{S}_{1}$ consists of the elements $gh\sim g^{2}h^{2}\sim
g^{3}h^{3}\sim g^{4}h^{4}\sim gh^{4}\sim g^{2}h^{3}\sim g^{3}h^{2}\sim
g^{4}h\,$ and $\mathcal{S}_{2}$ consists of $gh^{2}\sim g^{2}h^{4}\sim
g^{3}h\sim g^{4}h^{3}\sim gh^{3}\sim g^{2}h\sim g^{3}h^{4}\sim g^{4}h^{2}$
where $\mathcal{S}_{1}\cup \mathcal{S}_{2}
=\{g^{A}h^{B}|A,B\neq 0\ {\rm mod}\ 5\}$. The irrational GMFs occur when the elements of the two sets $\mathcal{S%
}_{1},\mathcal{S}_{2}$ are of different Fricke type, e.g. $\ \mathcal{S}_{1}$
is Fricke but $\mathcal{S}_{2}$ is non-Fricke. Indeed, otherwise we have
some of the cases, analysed previously in \cite{IT} which lead to rational
GMFs.

\begin{prop}
\label{order_h=5} For $g=5+$ and $h\in G_{g}\equiv \mathrm{HN}$ of order $5$%
, the following class structures give rise to a genus zero fixing group $%
\tilde{\Gamma}_{h,g}$ for irrational GMFs $Z\left[ 
\begin{array}{c}
h \\ 
g
\end{array}
\right] (5\tau )$ as follows:

\item[(i)]  $g$ Fricke, $h$ non-Fricke, $\mathcal{S}_{1}$ Fricke, $\mathcal{S%
}_{2}$ non-Fricke, 
\[
\tilde{\Gamma}_{h,g}=\;\langle \Gamma _{0}(25),T^{1/5}W_{25}\widetilde{%
\delta }_{5}T^{1/5}\rangle .
\]

\item[(ii)]  $g$ Fricke, $h$ non-Fricke, $\mathcal{S}_{1}$ non-Fricke, $%
\mathcal{S}_{2}$ Fricke, 
\[
\tilde{\Gamma}_{h,g}=\;\langle \Gamma _{0}(25),T^{2/5}W_{25}T^{2/5}\rangle 
\]
where $T^{r/5}=\left( 
\begin{array}{cc}
1 & r/5 \\ 
0 & 1
\end{array}
\right) $, $\widetilde{\delta }_{5}=\left( 
\begin{array}{cc}
2 & 1 \\ 
25 & 13
\end{array}
\right) $ and $W_{25}=$ $\left( 
\begin{array}{cc}
0 & -1 \\ 
25 & 0
\end{array}
\right) $, the Fricke involution.
\end{prop}

\textbf{Proof:} Firstly we will prove that $\tilde{\Gamma}_{h,g}$ contains $%
\Gamma _{0}(25)$. Due to Lemma \ref{delta-p lemma} it is sufficient to
demonstrate invariance with respect to $\delta _{5}$.

According to (\ref{Delta like inv}) we have $Z\left[ 
\begin{array}{c}
h \\ 
g
\end{array}
\right] =Z\left[ 
\begin{array}{c}
h^{d} \\ 
g^{2}
\end{array}
\right] $ for some $d$. Clearly $h^{d}\in G_{g}$ is of order $5$ and has
irrational characters. From inspection of the ATLAS \cite{CCNPW} we have
only two such classes with algebraically conjugate irreducible characters
distinguished by $\left( \frac{d}{5}\right) =\pm 1$. Therefore there are two
possibilities: either $Z\left[ 
\begin{array}{c}
h \\ 
g
\end{array}
\right] =Z\left[ 
\begin{array}{c}
h \\ 
g^{2}
\end{array}
\right] $ or $Z\left[ 
\begin{array}{c}
h^{3} \\ 
g^{2}
\end{array}
\right] $. The first possibility can be ruled out since it implies after an $%
ST^{-1\text{ }}$transformation that $Z\left[ 
\begin{array}{c}
g^{-1} \\ 
gh
\end{array}
\right] =Z\left[ 
\begin{array}{c}
g^{-2} \\ 
g^{2}h
\end{array}
\right] $. This is impossible since $gh\in \mathcal{S}_{1}$ and $g^{2}h\in 
\mathcal{S}_{2}$ and $\mathcal{S}_{1}$ and $\mathcal{S}_{2}$ are of
different Fricke/non-Fricke type for both cases (i) and (ii) above.
Therefore the second possibility remains which implies $\delta _{5}$
invariance and the result follows.

Case (i): From the analysis in \cite{IT}, where all rational cases with $%
o(g)=o(h)=5$ are analysed, it follows that $Z\left[ 
\begin{array}{c}
g^{a}h^{b} \\ 
g^{2}h^{3}
\end{array}
\right] $ cannot be rational for $(a,b)\neq (0,0)\ {\rm mod}\ 5$. If we choose $%
a$,$b$ such that $g^{a}h^{b}\in G_{g^{2}h^{3}}$ ($\phi
_{g^{2}h^{3}}(g^{a}h^{b})=1$), this GMF must be quadratically irrational.
Since the Harada-Norton group has only one pair of irrational classes of
order $5$ so that $Z\left[ 
\begin{array}{c}
g^{a}h^{b} \\ 
g^{2}h^{3}
\end{array}
\right] $ is either equal to $Z\left[ 
\begin{array}{c}
h \\ 
g
\end{array}
\right] $ or its algebraic conjugate. We may then further restrict the
choice of $a,b$ so that $Z\left[ 
\begin{array}{c}
h \\ 
g
\end{array}
\right] =Z\left[ 
\begin{array}{c}
g^{a}h^{b} \\ 
g^{2}h^{3}
\end{array}
\right] $.

Firstly we observe that $b\neq 0$ mod $5$ since $g^{a}h^{b}$ should
be non-Fricke since $h$ is. If $a=0$ mod $5$ then $Z\left[ 
\begin{array}{c}
h \\ 
g
\end{array}
\right] =Z\left[ 
\begin{array}{c}
h^{b} \\ 
g^{2}h^{3}
\end{array}
\right] $ and after a $T^{-1\text{ }}$transformation $Z\left[ 
\begin{array}{c}
gh \\ 
g
\end{array}
\right] =Z\left[ 
\begin{array}{c}
g^{2}h^{3+b} \\ 
g^{2}h^{3}
\end{array}
\right] $ and therefore $g^{2}h^{3+b}$ is also Fricke and hence $b\neq 1,3$ mod $5$. Applying a $T^{\text{ }}$transformation we can also see
that $g^{3}h^{2+b}$ is Fricke and hence $b\neq 2,4$ mod $5$. Thus $%
a\neq 0$ mod $5$.

Since $g^{a}h^{b}$ is non-Fricke we conclude that $(a,b)\notin \{(s,ns)|$ $%
\left( \frac{n}{5}\right) =1,$ $s\neq 0$ mod $5\}$. Applying a $%
T^{\pm 1\text{ }}$transformation we also conclude that $(a,b)\notin \{(\pm
2,s),(s\pm 2,\;ns\pm 3)|\left( \frac{n}{5}\right) =-1,s\neq 0$ mod 5.
Thus either $(a,b)=$ $(1,2)$ or $(4,3)$ both leading to the same GMF since $%
(gh^{2})^{-1}=g^{4}h^{3}$ and all characters are real.

Thus we have shown that $Z\left[ 
\begin{array}{c}
h \\ 
g
\end{array}
\right] (\tau )=Z\left[ 
\begin{array}{c}
gh^{2} \\ 
g^{2}h^{3}
\end{array}
\right] (\tau )=Z\left[ 
\begin{array}{c}
h \\ 
g
\end{array}
\right] (\alpha \tau )$ with $\alpha \equiv TS\delta _{5}T$ . Furthermore we
find that all the singular cusps of the GMF are identified under $\Gamma
_{h,g}=\langle \Gamma _{0}^{0}(5),\alpha \rangle $ which is therefore a
genus zero fixing group. Clearly $\Gamma _{h,g}$ is a subgroup of $\Gamma $
where $\alpha ^{3}=$ $\delta _{5}^{2}=(\alpha \delta _{5})^{2}=1\ {\rm mod}\ \Gamma (5)$ and so $\Gamma _{h,g}/\Gamma (5)=D_{3}$ which is a (maximal)
subgroup of $L_{2}(5)$. The corresponding hauptmodul for $\widetilde{\Gamma }%
_{h,g}=\langle \Gamma _{0}(25),\tilde{\alpha}\rangle $ can be explicitly
expressed as 
\begin{eqnarray}
Z\left[ 
\begin{array}{c}
h \\ 
g
\end{array}
\right] (5\tau ) &=&\frac{\eta (\tau )\eta (\tau +2/5)\eta (\tau +3/5)}{\eta
(\tau +1/5)\eta (\tau +4/5)\eta (25\tau )}+1-\sqrt{5}=  \nonumber \\
&=&q^{-1}+0+(\frac{3}{2}-\frac{5\sqrt{5}}{2})q-10q^{2}+5q^{3}  \nonumber \\
&&+(21+5\sqrt{5})q^{4}+(\frac{-25}{2}+\frac{25\sqrt{5}}{2})q^{5}+...
\label{HN order 5 }
\end{eqnarray}
which is known as the $25\symbol{126}b$ series \cite{N3}.

Assuming that the GMF is replicable and based on numerical matching L. Queen
conjectures that the head characters for $h\in G_{5+}\equiv \mathrm{HN}$ can
be expanded in terms of the irreducible characters $\chi _{i}(h)$ of HN (in
ATLAS notation \cite{CCNPW}) to give \cite{Q1}, \cite{Q2} 
\begin{eqnarray}
Z\left[ 
\begin{array}{c}
h \\ 
g
\end{array}
\right] (5\tau ) &=&\frac{1}{q}+0+(\chi _{1}(h)+\chi _{3}(h))q+\chi
_{4}(h)q^{2}+(\chi _{1}(h)+\chi _{5}(h))q^{3}+  \nonumber \\
&&(\chi _{1}(h)+\chi _{2}(h)+\chi _{5}(h)+\chi _{6}(h))q^{4}+  \nonumber \\
&&(\chi _{1}(h)+\chi _{2}(h)+\chi _{4}(h)+\chi _{5}(h)+\chi
_{11}(h))q^{5}+...  \label{HN head ch exp}
\end{eqnarray}
and therefore (\ref{HN order 5 }) corresponds to the head character
expansion (\ref{HN head ch exp}) for the $5C$ class of $\mathrm{HN}$ \cite
{CCNPW}.

Case (ii). This class structure is obtained from that in Case (i) by
replacing $h$ by $h^{2}$. Therefore the fixing group given by (\ref{Gamma*k
f.gr.}) with $k=2$ is $\tilde{\Gamma}_{h,g}=\langle \Gamma
_{0}(25),T^{2/5}W_{25}T^{2/5}\rangle $. Its hauptmodul is then $\frac{\eta
(\tau )\eta (\tau +1/5)\eta (\tau +4/5)}{\eta (\tau +2/5)\eta (\tau
+3/5)\eta (25\tau )}+1+\sqrt{5}$. This is known as the $25\symbol{126}a$
series \cite{N3} and has algebraically conjugate coefficients to those of
case (i). The head character expansion corresponds to the $5D$ class of $%
\mathrm{HN}$ (\ref{HN head ch exp}).$\qed $

\subsubsection{Irrational GMFs for $g=7+$ and $o(h)=7$}

The Held group has two pairs of quadratically irrational classes of order $7$%
. In (\ref{representatives}) $m=3$, $N=2$. The two subsets of powers for
each $h$ of order $7$ (\ref{h-a sets}), conjugate in $G_{7+}=\mathrm{He}$
are: $\left\{ h^{a}\right\} $ for $\left( \frac{a}{7}\right) =1$ and $%
\left\{ h^{a}\right\} $ for $\left( \frac{a}{7}\right) =-1$ \cite{CCNPW}.
Therefore $gh\stackrel{\bf M}{\sim }gh^{2}\stackrel{\bf M}{\sim }gh^{4}$
and $gh^{3}\stackrel{\bf M}{\sim }gh^{5}\stackrel{\bf M}{\sim }gh^{6}$.
Since $g_{0}\stackrel{\bf M}{\sim }g_{0}^{s}$ ($s=1$,$2$,...,$p-1$) for
any Monster element $g_{0}$, we have the following two disjoint sets $%
\mathcal{S}_{1},\mathcal{S}_{2}$ of conjugate elements in $\bf{M}$ defined
by

\begin{eqnarray*}
\mathcal{S}_{1} &\equiv &\{(gh^{n})^{s}|\left( \frac{n}{7}\right) =1,s\neq 0\ {\rm mod}\ 7\} \\
\mathcal{S}_{2} &\equiv &\{(gh^{n})^{s}|\left( \frac{n}{7}\right) =-1,s\neq 0\ {\rm mod}\ 7\}
\end{eqnarray*}
where $\mathcal{S}_{1}\cup \mathcal{S}_{2}=\{g^{A}h^{B}|A,B\neq 0\ {\rm mod}\ 7\}$.

\begin{prop}
\label{order 7}For $g=7+$ and $h$ of order $7$, $h\in G_{7+}\equiv \mathrm{He%
}$ the following class structures give rise to a genus zero fixing group $%
\tilde{\Gamma}_{h,g}$ for irrational GMFs $Z\left[ 
\begin{array}{c}
h \\ 
g
\end{array}
\right] (7\tau )$:

\item[(i)]  $g$ Fricke, $h$ non-Fricke, the set $\mathcal{S}_{1}$ Fricke, $%
\mathcal{S}_{2}$ non-Fricke, 
\[
\tilde{\Gamma}_{h,g}=\;\langle \Gamma _{0}(49),\widetilde{\alpha }\rangle 
\text{, }\widetilde{\alpha }\equiv \left( 
\begin{array}{cc}
3 & 2/7 \\ 
-35 & -3
\end{array}
\right) 
\]

\item[(ii)]  $g$ Fricke, $h$ non-Fricke, $\mathcal{S}_{1}$ non-Fricke, $%
\mathcal{S}_{2}$ Fricke, 
\[
\tilde{\Gamma}_{h,g}=\;\langle \Gamma _{0}(49),\widetilde{\alpha }%
^{*}\rangle \text{, }\widetilde{\alpha }^{*}\equiv \left( 
\begin{array}{cc}
3 & -2/7 \\ 
35 & -3
\end{array}
\right) 
\]

\item[(iii)]  $g$ Fricke, $h$ Fricke, $\mathcal{S}_{1}$ and $\mathcal{S}_{2}$
sets Fricke. There are two possible fixing groups, which are isomorphic and
usually both denoted by $7||7+$.
\end{prop}

\textbf{Proof}: Firstly we will prove that $\tilde{\Gamma}_{h,g}$ contains $%
\Gamma _{0}(49)$. Due to Lemma \ref{delta-p lemma} it is sufficient to
demonstrate invariance with respect to $\delta _{7}$.

According to (\ref{Delta like inv}) we have $Z\left[ 
\begin{array}{c}
h \\ 
g
\end{array}
\right] =Z\left[ 
\begin{array}{c}
h^{d} \\ 
g^{2}
\end{array}
\right] $ for some $d$. We need to show that $Z\left[ 
\begin{array}{c}
h \\ 
g
\end{array}
\right] (\tau )=Z\left[ 
\begin{array}{c}
h^{4} \\ 
g^{2}
\end{array}
\right] (\tau )\equiv Z\left[ 
\begin{array}{c}
h \\ 
g
\end{array}
\right] (\delta _{7}\tau )$. Assume that $\left( \frac{d}{7}\right) =-1$.
Then we have $Z\left[ 
\begin{array}{c}
h \\ 
g
\end{array}
\right] =Z\left[ 
\begin{array}{c}
h^{-1} \\ 
g^{2}
\end{array}
\right] =Z\left[ 
\begin{array}{c}
h \\ 
g^{4}
\end{array}
\right] =Z\left[ 
\begin{array}{c}
h^{-1} \\ 
g^{8}
\end{array}
\right] \equiv Z\left[ 
\begin{array}{c}
h^{-1} \\ 
g
\end{array}
\right] $ which is impossible since $h$ and $h^{-1}$ are not conjugate in $%
G_{g}$. Therefore $\left( \frac{d}{7}\right) =1$ and the statement follows.

When $p=7$ there are two pairs of conjugate irrational characters in $%
G_{7+}\equiv \mathrm{He}$ \cite{CCNPW}. If we consider a GMF on another
Fricke twisted sector, for convenience, a $g^{3}h^{5}$ twisted sector, then
for a pair $(a,b)\neq (0,0)$ mod $7$ and $g^{a}h^{b}\in
G_{g^{3}h^{5}}$ ($\phi _{g^{3}h^{5}}(g^{a}h^{b})=1$), $Z\left[ 
\begin{array}{c}
g^{a}h^{b} \\ 
g^{3}h^{5}
\end{array}
\right] $ must be quadratically irrational and $\langle g,h\rangle =\langle
g^{3}h^{5},g^{a}h^{b}\rangle $, so that it is either equal to $Z\left[ 
\begin{array}{c}
h \\ 
g
\end{array}
\right] $ or its algebraic conjugate. By choosing an appropriate multiple of 
$(a,b)$ we may further restrict the choice of $a,b$ so that 
\begin{equation}
Z\left[ 
\begin{array}{c}
h \\ 
g
\end{array}
\right] =Z\left[ 
\begin{array}{c}
g^{a}h^{b} \\ 
g^{3}h^{5}
\end{array}
\right]  \label{Z g3h5}
\end{equation}

Case (i) Suppose $a=0$ mod $7$. Then $T^{-1}$ transformation gives $%
Z\left[ 
\begin{array}{c}
gh \\ 
g
\end{array}
\right] =Z\left[ 
\begin{array}{c}
g^{3}h^{5+b} \\ 
g^{3}h^{5}
\end{array}
\right] $ and since $gh$ is Fricke, so is $g^{3}h^{5+b}$ i.e. it is in the $%
\mathcal{S}_{1}$ Fricke set and therefore $b\neq 4,3,6$ mod  $7$.
Since $h^{4}\stackrel{}{\sim }h^{2}\stackrel{}{\sim }h$ and $h^{3}\stackrel{%
}{\sim }h^{5} $ in the corresponding centraliser, 
$b\neq 1,2,5$ mod  $7$. Thus $a\neq 0$ mod  $7$.

Since $g^{a}h^{b}\in \mathcal{S}_{2}$ is non-Fricke, $b\neq 0$ mod  $%
7$, $(a,b)\notin \{(s,ns)|$ $\left( \frac{n}{7}\right) =1$, 
$s\neq 0$ mod $7\}$. 
Furthermore, a $T^{-1}$ transformation of (\ref{Z g3h5}) gives $%
Z\left[ 
\begin{array}{c}
gh \\ 
g
\end{array}
\right] =Z\left[ 
\begin{array}{c}
g^{3+a}h^{5+b} \\ 
g^{3}h^{5}
\end{array}
\right] $, and since $gh$ is Fricke, so is $g^{3+a}h^{5+b}$ thus $%
(a,b)\notin \{(s-3,ns-5)$, $(4,s)|$ $\left( \frac{n}{7}\right) =-1,$ $s\neq 0 $ mod  $7\}$. If a given pair $(a,b)$ is inadmissible, (i.e. does
not satisfy (\ref{Z g3h5})) so are the pairs leading to conjugates in the
same set: $(na,nb)$, $\left( \frac{n}{7}\right) =1$. Therefore $(a,b)=$ $%
(5n,4n)$, $\left( \frac{n}{7}\right) =1$. Hence $Z\left[ 
\begin{array}{c}
h \\ 
g
\end{array}
\right] (\tau )=Z\left[ 
\begin{array}{c}
g^{5}h^{4} \\ 
g^{3}h^{5}
\end{array}
\right] (\tau )=Z\left[ 
\begin{array}{c}
h \\ 
g
\end{array}
\right] (\alpha \tau )$ where $\alpha =\left( 
\begin{array}{cc}
3 & 2 \\ 
-5 & -3
\end{array}
\right) $ and thus the GMF is $\alpha $ modular invariant. Furthermore we
find that all the singular cusps of the GMF are identified under $\Gamma
_{h,g}=\langle \Gamma _{0}^{0}(7),\alpha \rangle $ which is therefore a
genus zero fixing group. Introducing $\beta =\delta _{7}\alpha ^{-1}=\left( 
\begin{array}{cc}
1 & 3 \\ 
-1 & -2
\end{array}
\right) $, one can check that $\beta ^{3}=$ $\alpha ^{2}=(\beta \alpha
)^{3}=1\ {\rm mod} \ \Gamma (7)$. Therefore $\langle \Gamma (7),\delta
_{7},\alpha \rangle /\Gamma (7)\simeq A_{4}$, a subgroup of $L_{2}(7)$. The
hauptmodul for the genus zero fixing group $\tilde{\Gamma}_{h,g}=\langle
\Gamma _{0}(49),\tilde{\alpha}\rangle $ is 
\begin{eqnarray}
Z\left[ 
\begin{array}{c}
h \\ 
g
\end{array}
\right] (7\tau ) &=&\frac{\eta (\tau )\eta (\tau +3/7)\eta (\tau +5/7)\eta
(\tau +6/7)}{\eta ^{4}(7\tau )}+\frac{1}{2}-i\frac{\sqrt{7}}{2}=  \nonumber
\\
&&\frac{1}{q}+0+(-\frac{3}{2}+i\frac{\sqrt{7}}{2})q+(-\frac{5}{2}-i\frac{3%
\sqrt{7}}{2})q^{2}+2q^{3}+  \nonumber \\
&&(3-i\sqrt{7})q^{4}-3q^{5}+...  \label{He hauptmodul}
\end{eqnarray}
known as the $49\symbol{126}a$ series \cite{N3}.

Assuming that the GMF is replicable and based on numerical matching L. Queen
conjectures that the head characters for $h\in G_{7+}\equiv \mathrm{He}$ can
be expanded in terms of the irreducible characters $\chi _{i}(h)$ of $%
\mathrm{He}$ (in ATLAS notation \cite{CCNPW}) to give \cite{Q1},\cite{Q2} 
\begin{eqnarray}
Z\left[ 
\begin{array}{c}
h \\ 
g
\end{array}
\right] (7\tau ) &=&\frac{1}{q}+0+\chi _{2}(h)q+(\chi _{3}(h)+\chi
_{4}(h))q^{2}+(\chi _{1}(h)+\chi _{6}(h))q^{3}+  \nonumber \\
&&(\chi _{1}(h)+\chi _{6}(h)+\chi _{11}(h))q^{4}+  \nonumber \\
&&(\chi _{1}(h)+\chi _{2}(h)+\chi _{3}(h)+\chi _{6}(h)+\chi
_{14}(h))q^{5}+...  \label{He char exp}
\end{eqnarray}
Thus the expansion (\ref{He hauptmodul}) is associated with the $7E$ class
of $\mathrm{He.}$

Case (ii). This class structure is obtained from that in Case (i) by
replacing $h$ by $h^{-1}$. Therefore the fixing group is given by (\ref
{Gamma*k f.gr.}) with $k=-1$ mod  $7$. The fixing group contains $%
\alpha ^{*}=\left( 
\begin{array}{cc}
3 & -2 \\ 
5 & -3
\end{array}
\right) $ i.e. $\Gamma _{h,g}=\langle \Gamma _{0}^{0}(7),\alpha ^{*}\rangle $%
. $h$ is in $7D$ class of $\mathrm{He}$, with algebraically conjugate
characters to those of $7E$. The GMF $Z\left[ 
\begin{array}{c}
h \\ 
g
\end{array}
\right] (7\tau )=\frac{\eta (\tau )\eta (\tau +1/7)\eta (\tau +2/7)\eta
(\tau +4/7)}{\eta ^{4}(7\tau )}$ $+\frac{1}{2}+i\frac{\sqrt{7}}{2}$
corresponding to the $49\symbol{126}b$ series \cite{N3}, is a hauptmodul for
the genus zero fixing group $\tilde{\Gamma}_{h,g}=\langle \Gamma _{0}(49),%
\tilde{\alpha}^{*}\rangle $. This hauptmodul has algebraically conjugate
coefficients compared to (\ref{He hauptmodul}).

Case (iii) From earlier remarks there are two algebraically conjugate order $%
7$ classes of $\mathrm{He}$ remaining. We will demonstrate $S$-invariance.
Since $Z\left[ 
\begin{array}{c}
h \\ 
g
\end{array}
\right] =Z\left[ 
\begin{array}{c}
h \\ 
g^{2}
\end{array}
\right] $, after an $S$-transformation $Z\left[ 
\begin{array}{c}
g^{-1} \\ 
h
\end{array}
\right] =Z\left[ 
\begin{array}{c}
g^{-2} \\ 
h
\end{array}
\right] $ and thus $\phi _{h}(g^{-1})=\phi _{h}(g^{-2})$ i.e. $\phi
_{h}(g)=1 $. So $g\in G_{h}\simeq G_{g}$. We have only one pair of
irrational classes left, hence $Z\left[ 
\begin{array}{c}
g \\ 
h
\end{array}
\right] =Z\left[ 
\begin{array}{c}
h \\ 
g
\end{array}
\right] $, or $Z\left[ 
\begin{array}{c}
g \\ 
h
\end{array}
\right] =Z\left[ 
\begin{array}{c}
h^{-1} \\ 
g
\end{array}
\right] $, its conjugate. Assume that $Z\left[ 
\begin{array}{c}
g \\ 
h
\end{array}
\right] =Z\left[ 
\begin{array}{c}
h \\ 
g
\end{array}
\right] $, then after an $ST^{-1}$ transformation $Z\left[ 
\begin{array}{c}
g^{-1} \\ 
gh
\end{array}
\right] =Z\left[ 
\begin{array}{c}
h^{-1} \\ 
gh
\end{array}
\right] $, which implies that $Z\left[ 
\begin{array}{c}
g^{A} \\ 
gh
\end{array}
\right] =Z\left[ 
\begin{array}{c}
h^{A} \\ 
gh
\end{array}
\right] $ for any $A$. Taking $A=-2$, after a $T^{-1}$ transformation we
obtain $Z\left[ 
\begin{array}{c}
g^{-1}h \\ 
gh
\end{array}
\right] =Z\left[ 
\begin{array}{c}
(g^{-1}h)^{-1} \\ 
gh
\end{array}
\right] $. This is impossible since the irrational characters of $\mathrm{He}$
are complex. Hence $Z\left[ 
\begin{array}{c}
h \\ 
g
\end{array}
\right] (\tau )=Z\left[ 
\begin{array}{c}
g^{-1} \\ 
h
\end{array}
\right] (\tau )=Z\left[ 
\begin{array}{c}
h \\ 
g
\end{array}
\right] (S\tau )$.

Consider the GMF $Z^{*}$, algebraically conjugate to $Z$%
\begin{equation}
Z^{*}\left[ 
\begin{array}{c}
h \\ 
g
\end{array}
\right] \equiv Z\left[ 
\begin{array}{c}
h^{-1} \\ 
g
\end{array}
\right] .  \label{Z star}
\end{equation}
From (\ref{Z g3h5}) and $S$ invariance we have $Z\left[ 
\begin{array}{c}
h^{-1} \\ 
g
\end{array}
\right] =$ $Z\left[ 
\begin{array}{c}
h^{-2} \\ 
g
\end{array}
\right] =$ $Z\left[ 
\begin{array}{c}
g \\ 
h^{2}
\end{array}
\right] =$ $Z\left[ 
\begin{array}{c}
g^{b}h^{2a} \\ 
g^{5}h^{6}
\end{array}
\right] =$ $Z\left[ 
\begin{array}{c}
g^{b}h^{2a} \\ 
(g^{5}h^{6})^{2}
\end{array}
\right] =$ $\;Z\left[ 
\begin{array}{c}
g^{b}h^{2a} \\ 
g^{3}h^{5}
\end{array}
\right] =$ $Z^{*}\left[ 
\begin{array}{c}
g^{-b}h^{-2a} \\ 
g^{3}h^{5}
\end{array}
\right] $ i.e. if $Z$ satisfies (\ref{Z g3h5}) for a pair $(a,b)$ then $%
Z^{*} $ satisfies (\ref{Z g3h5}) for the pair $(-b,-2a)$ (and therefore $%
(-nb,-2na) $ with $\left( \frac{n}{7}\right) =1$).

We can show that $Z$ and $Z^{*}$ cannot satify (\ref{Z g3h5}) for one and
the same pair $(a,b)$ (otherwise they would have the same fixing group and
possibly hauptmodul). Let us assume that it is possible and therefore $%
(a,b)=(-2b,-4a)$ mod  $7$ ($n=2$): i.e. $(a,b)=$ $(s,3s)$, $s\neq 0$ 
mod  $7$. For $s=1$ we have $Z\left[ 
\begin{array}{c}
h \\ 
g
\end{array}
\right] =Z\left[ 
\begin{array}{c}
gh^{3} \\ 
g^{3}h^{5}
\end{array}
\right] $. Since $Z\left[ 
\begin{array}{c}
h \\ 
g
\end{array}
\right] =Z\left[ 
\begin{array}{c}
h^{3} \\ 
g^{5}
\end{array}
\right] $ it follows that $Z\left[ 
\begin{array}{c}
h \\ 
g
\end{array}
\right] =Z\left[ 
\begin{array}{c}
g^{5}h^{2} \\ 
gh
\end{array}
\right] $. Applying two times the last transformation we get $Z\left[ 
\begin{array}{c}
h \\ 
g
\end{array}
\right] =Z\left[ 
\begin{array}{c}
g^{4}h^{5} \\ 
h^{5}
\end{array}
\right] $ which is a contradiction since clearly $\phi
_{h^{5}}(g^{4}h^{5})\neq 1$. Thus $(a,b)\neq $ $(s,3s)$ for $\left( \frac{s}{%
7}\right) =1$. For $s=3$ we have $Z\left[ 
\begin{array}{c}
h \\ 
g
\end{array}
\right] =Z\left[ 
\begin{array}{c}
g^{3}h^{2} \\ 
g^{3}h^{5}
\end{array}
\right] $. Applying two times this transformation we get $Z\left[ 
\begin{array}{c}
h \\ 
g
\end{array}
\right] =Z\left[ 
\begin{array}{c}
g^{4}h \\ 
h^{2}
\end{array}
\right] $ which is a contradiction since clearly $\phi _{h^{2}}(g^{4}h)\neq
1 $. Thus $(a,b)\neq $ $(s,3s)$ for $\left( \frac{s}{7}\right) =-1$ and
finally $(a,b)\neq $ $(s,3s)$ in general.

Suppose $a=0$ mod  $7$. Then $Z\left[ 
\begin{array}{c}
g \\ 
h
\end{array}
\right] =Z\left[ 
\begin{array}{c}
g^{3}h^{5} \\ 
h^{b}
\end{array}
\right] $ which is impossible, because $\phi _{h}(g)=1$, $\phi _{h^{b}}(g)=1 
$, $\phi _{h^{b}}(h^{b})=\omega _{7}$ . Thus $a\neq 0$ mod  $7$ and
similarly $b\neq 0$ mod  $7$.

Since $\phi _{g^{3}h^{5}}(g^{3}h^{5})=\omega _{7}$, $g^{a}h^{b}$ must not be
a power of $g^{3}h^{5}$, $(a,b)\neq $ $(3s,5s)$. Applying once again (\ref{Z
g3h5}) we have 
\[
Z\left[ 
\begin{array}{c}
h \\ 
g
\end{array}
\right] =Z\left[ 
\begin{array}{c}
g^{a(3+b)}h^{5a+b^{2}} \\ 
g^{2-2a}h^{1-2b}
\end{array}
\right] . 
\]
Let $a=4b^{2}$ mod  $7$ and $b=\pm 1$ or $\pm 2$ mod  $7$.
Then $Z\left[ 
\begin{array}{c}
h \\ 
g
\end{array}
\right] =Z\left[ 
\begin{array}{c}
g^{5b^{2}+4b^{3}} \\ 
g^{2-b^{2}}h^{1-2b}
\end{array}
\right] $ is possible only if $\phi _{g^{2-b^{2}}h^{1-2b}}(g)=1$. However
according to Theorem \ref{theorem_proj1}, since $\phi _{h}(g)=1$ it is not
possible to have $\phi _{g^{2-b^{2}}h^{1-2b}}(g)=1$ unless $2-b^{2}=0$ mod $7$ which is not the case. Therefore $(a,b)\neq $ $(4,\pm 1)$, $%
(2,\pm 2)$, and we exclude also their conjugates of the form $(na,nb)$, $%
\left( \frac{n}{7}\right) =1$: $(a,b)\notin \{(4n,\pm n)$, $(2n,\pm 2n)|$ $%
\left( \frac{n}{7}\right) =1\}$. Finally we exclude all pairs, which can be
represented as $(-b,-2a)$ from an already excluded pair $(a,b)$. It is due
to the fact that if $(a,b)$ satisfies (\ref{Z g3h5}), $(-b,-2a)$ must
satisfy the same relation for $Z^{*}$ (\ref{Z star}). The remaining pairs
are $(a,b)=(5n,4n)$, $\left( \frac{n}{7}\right) =1$ ($(5,4)$, $(3,1)$ and $%
(6,2)$); their $Z^{*}$-partners of the form $(-b,-2a)$: $(3,4)$, $(6,1)$, $%
(5,2)$. This leads to modular invariance under $\alpha =\left( 
\begin{array}{cc}
3 & 2 \\ 
-5 & -3
\end{array}
\right) $ for $Z$ (for $(a,b)=(5,4)$), and $\alpha ^{\prime }=\left( 
\begin{array}{cc}
-4 & 1 \\ 
-5 & 1
\end{array}
\right) $ for $Z^{*}$ (for $(-b,-2a)=(3,4)$). Note that $\alpha $ is of
order $2$ in $\Gamma (7)$ and $\alpha ^{\prime }$ is of order $4$ in $\Gamma
(7)$ i.e. the mapping to $g^{3}h^{5}$ Fricke twisted sector (\ref{Z g3h5})
requires modular transformations of different order in each case. Each of
the two fixing groups $\Gamma _{h,g}=\langle \Gamma _{0}^{0}(7),\alpha
,S\rangle $, $\langle \Gamma _{0}^{0}(7),\alpha ^{\prime },S\rangle $
matches all the singular cusps of the GMF and therefore each group is of
genus zero. These two fixing groups are isomorphic. Indeed, introducing $%
\alpha _{1}=\alpha S$, $\alpha _{2}=\delta _{7}S$, $\alpha _{3}=S$, we have $%
\alpha _{1}^{2}=\alpha _{2}^{2}=\alpha _{3}^{2}=$ $(\alpha _{1}\alpha
_{2})^{3}=(\alpha _{2}\alpha _{3})^{3}=(\alpha _{3}\alpha _{1})^{2}=$ $1\ 
{\rm mod} \ \Gamma (7)$, the defining relations of the group $S_{4}$.
Introducing $\beta _{1}=S$, $\beta _{2}=S\alpha ^{\prime }$, we have $\beta
_{1}^{2}=\beta _{2}^{3}=(\beta _{1}\beta _{2})^{4}=1\ {\rm mod} \ \Gamma (7)$%
, also defining the group $S_{4}$, a maximal subgroup of $L_{2}(7)$.
Therefore $\langle \Gamma (7),\delta _{7},\alpha ,S\rangle /\Gamma (7)\simeq
S_{4}$ and $\langle \Gamma (7),\alpha ^{\prime },S\rangle /\Gamma (7)\simeq
S_{4}$. According to the definition in \cite{FMN} for both groups we use the
same notation $\tilde{\Gamma}_{h,g}=7||7+$. The GMFs represents the
hauptmoduls for these groups and can be found from those of Case (i) with
symmetrization: 
\begin{eqnarray*}
Z\left[ 
\begin{array}{c}
h \\ 
g
\end{array}
\right] (7\tau ) &=&\frac{\eta (\tau )\eta (\tau +3/7)\eta (\tau +5/7)\eta
(\tau +6/7)}{\eta ^{4}(7\tau )}+ \\
&&i\sqrt{7}\frac{\eta (49\tau )\eta (\tau +1/7)\eta (\tau +2/7)\eta (\tau
+4/7)}{\eta ^{4}(7\tau )} \\
&=&q^{-1}+0+(-\frac{3}{2}+i\frac{3\sqrt{7}}{2})q+(1-i\sqrt{7})q^{2}+9q^{3}+
\\
&&3(1-i\sqrt{7})q^{4}+4q^{5}+...
\end{eqnarray*}

\begin{eqnarray*}
Z^{*}\left[ 
\begin{array}{c}
h \\ 
g
\end{array}
\right] (7\tau ) &=&\frac{\eta (\tau )\eta (\tau +1/7)\eta (\tau +2/7)\eta
(\tau +4/7)}{\eta ^{4}(7\tau )}+ \\
&&i\sqrt{7}\frac{\eta (49\tau )\eta (\tau +3/7)\eta (\tau +5/7)\eta (\tau
+6/7)}{\eta ^{4}(7\tau )}
\end{eqnarray*}

$Z^{*}$ has a series expansion with complex conjugate coefficients to $Z$.
According to (\ref{He char exp}) $Z$ is the GMF when $h$ is an element of $%
7A $ class of $\mathrm{He}$, $Z^{*}$ - when $h$ is an element of $7B$ class
of $\mathrm{He}$.$\qed $

\subsubsection{Irrational GMFs for $g=11+$ and $o(h)=11$}

The \textrm{M}$_{12}$ group has one pair of quadratically irrational classes
of order $11$. In (\ref{representatives}) $m=2$, $N=2$. There are two
subsets of powers of $h$ (\ref{h-a sets}) conjugate in $G_{11+}\equiv $%
\textrm{M}$_{12}$: $\left\{ h^{a}\right\} $ for $\left( \frac{a}{11}\right)
=1$ and $\left\{ h^{a}\right\} $ for $\left( \frac{a}{11}\right) =-1$ \cite
{CCNPW}.

\begin{prop}
\label{h order 11}For $g=11+$, $h\in G_{11+}\equiv \mathrm{M}_{12}$ and $%
o(h)=11$, the class structure where all elements of the group $\langle
g,h\rangle $ are Fricke gives rise to a pair of irrational GMFs with
isomorphic genus zero fixing groups (usually both denoted by $11||11+$).
\end{prop}

\textbf{Proof}: Firstly we will prove that $\tilde{\Gamma}_{h,g}$ contains $%
\Gamma _{0}(121)$. Due to Lemma \ref{delta-p lemma} it is sufficient to
demonstrate invariance with respect to $\delta _{11}$.

According to (\ref{Delta like inv}) $Z\left[ 
\begin{array}{c}
h \\ 
g
\end{array}
\right] =Z\left[ 
\begin{array}{c}
h^{d} \\ 
g^{4}
\end{array}
\right] $ for some $d$. We need to show that $Z\left[ 
\begin{array}{c}
h \\ 
g
\end{array}
\right] (\tau )=Z\left[ 
\begin{array}{c}
h^{3} \\ 
g^{4}
\end{array}
\right] (\tau )\equiv Z\left[ 
\begin{array}{c}
h \\ 
g
\end{array}
\right] (\delta _{11}\tau )$. Clearly $h^{d}\in G_{g}$ is of order $11$ and
has irrational characters. From inspection of the ATLAS \cite{CCNPW} we have
only two such classes with algebraically conjugate irreducible characters
distinguished by $\left( \frac{d}{11}\right) =\pm 1$. Assume that $\left( 
\frac{d}{11}\right) =-1$. Then we have $Z\left[ 
\begin{array}{c}
h \\ 
g
\end{array}
\right] =Z\left[ 
\begin{array}{c}
h^{-1} \\ 
g^{4}
\end{array}
\right] =Z\left[ 
\begin{array}{c}
h^{(-1)^{5}} \\ 
g^{4^{5}}
\end{array}
\right] =Z\left[ 
\begin{array}{c}
h^{-1} \\ 
g
\end{array}
\right] $ which is impossible since $h$ and $h^{-1}$ are not conjugate in $%
G_{g}$. Therefore $\left( \frac{d}{11}\right) =1$, and the statement
follows. One can demonstrate also that $\phi _{h}(g)=1$ and hence $S$
invariance follows as in Proposition \ref{order 7}(iii).

Let us now consider a GMF related to another Fricke twisted sector, for
convenience, say a $gh^{2}$ twisted sector. When $o(h)=11$ there is only one
pair of elements in $\mathrm{M}_{12}$ with irrational characters \cite{CCNPW}%
. So $Z\left[ 
\begin{array}{c}
g^{a}h^{b} \\ 
gh^{2}
\end{array}
\right] $ cannot be rational for $(a,b)\neq (0,0)$ mod  $11.$ If we
choose $a$,$b$ such that $g^{a}h^{b}\in G_{gh^{2}}$ ($\phi
_{gh^{2}}(g^{a}h^{b})=1$), this GMF must be quadratically irrational. Since $%
\mathrm{M}_{12}$ has only one pair of irrational classes of order $11$ (and
no rational ones), then $Z\left[ 
\begin{array}{c}
g^{a}h^{b} \\ 
gh^{2}
\end{array}
\right] $ is either equal to $Z\left[ 
\begin{array}{c}
h \\ 
g
\end{array}
\right] $ or its algebraic conjugate. We may then further restrict the
choice of $a,b$ so that 
\begin{equation}
Z\left[ 
\begin{array}{c}
h \\ 
g
\end{array}
\right] =Z\left[ 
\begin{array}{c}
g^{a}h^{b} \\ 
gh^{2}
\end{array}
\right] .  \label{Z-gh2}
\end{equation}

Using (\ref{Z-gh2}) and $S$ invariance we have $Z\left[ 
\begin{array}{c}
h^{-1} \\ 
g
\end{array}
\right] =Z\left[ 
\begin{array}{c}
h^{-4} \\ 
g
\end{array}
\right] =$ $Z\left[ 
\begin{array}{c}
h^{-4.(-5)} \\ 
g^{-5}
\end{array}
\right] =$ $Z\left[ 
\begin{array}{c}
g^{-5} \\ 
h^{2}
\end{array}
\right] =$ $Z\left[ 
\begin{array}{c}
g^{-5b}h^{2a} \\ 
h^{2}g^{-5.2}
\end{array}
\right] =$ $Z\left[ 
\begin{array}{c}
(g^{5b}h^{-2a})^{-1} \\ 
gh^{2}
\end{array}
\right] $ and hence 
\[
Z^{*}\left[ 
\begin{array}{c}
h \\ 
g
\end{array}
\right] \equiv Z^{*}\left[ 
\begin{array}{c}
g^{5b}h^{-2a} \\ 
gh^{2}
\end{array}
\right] 
\]
where the definition of $Z^{*}$ is as in (\ref{Z star}). Thus if $(a,b)$ is
an invariance for $Z$ as in (\ref{Z-gh2}) then $(5b,-2a)$ is an invariance
for $Z^{*}$. Using similar techniques as in Proposition \ref{order 7}(iii)
we can prove that $(a,b)\neq (5b,-2a)$ mod  $11$ i.e. $(a,b)\neq $ $%
(s,9s)$.
Furthermore, since $\phi _{gh^{2}}(gh^{2})=\omega _{11}$, $(a,b)\neq $ $%
(s,2s)$.

As in Proposition \ref{order 7}(iii) one can show that $a,b\neq 0$ mod $11$. 
Applying once again the transformation (\ref{Z-gh2}) we have $%
Z\left[ 
\begin{array}{c}
h \\ 
g
\end{array}
\right] =Z\left[ 
\begin{array}{c}
g^{a(1+b)}h^{2a+b^{2}} \\ 
g^{2a+1}h^{2b+2}
\end{array}
\right] $. Let $a=5b^{2}$ mod  $11$ with $b^{2}\neq 1\ {\rm mod} $ $11$%
. Then $Z\left[ 
\begin{array}{c}
h \\ 
g
\end{array}
\right] =Z\left[ 
\begin{array}{c}
g^{5b^{2}(1+b)} \\ 
g^{1-b^{2}}h^{2b+2}
\end{array}
\right] $ is possible only if $\phi _{g^{1-b^{2}}h^{2b+2}}(g)=1$. However
according to theorem (\ref{theorem_proj1}), since $\phi _{h}(g)=1$ it is not
possible to have $\phi _{g^{1-b^{2}}h^{2b+2}}(g)=1$ unless $b^{2}=1$ mod $11$ 
which is not the case. Therefore $(a,b)\neq $ $(5b^{2},b)$ for $%
b=\pm 2$, $\pm 3$, $\pm 4$, $\pm 5$ mod $11$, and we exclude also
their conjugates of the form $(5b^{2}n,bn)$, $\left( \frac{n}{11}\right) =1$%
, $b\neq 0$, $\pm 1$ mod $11$. Finally we exclude all pairs which
can be represented as $(5b,-2a)$ from an already excluded pair $(a,b)$. It
is due to the fact that if $(a,b)$ satisfies (\ref{Z-gh2}), $(5b,-2a)$ must
satisfy the same relation for $Z^{*}$.

The remaining pairs are $(a,b)=(10n,10n)$, $\left( \frac{n}{11}\right) =1$
and their $Z^{*}$-partners of the form $(5b,-2a)$: $(6n,2n)$, $\left( \frac{n%
}{11}\right) =1$. This leads to modular invariance under $\alpha =\left( 
\begin{array}{cc}
1 & 1 \\ 
-2 & -1
\end{array}
\right) $ for $(a,b)=(10,10)$ -say for $Z$, and $\alpha ^{\prime }=\left( 
\begin{array}{cc}
1 & 5 \\ 
-2 & -9
\end{array}
\right) $ for $(a,b)=(6,2)$ for $Z^{*}$. The two fixing groups $\Gamma
_{h,g}=\langle \Gamma _{0}^{0}(11),\alpha ,S\rangle $, $\langle \Gamma
_{0}^{0}(11),\alpha ^{\prime },S\rangle $ are conjugate and isomorphic to $%
A_{5}$ ${\rm mod}\ \Gamma (11)$ (generated by $\alpha _{1}=\delta
_{11}S\alpha $, $\alpha _{2}=S$, $\alpha _{3}=\alpha $, which satisfy $%
\alpha _{1}^{2}=\alpha _{2}^{2}=\alpha _{3}^{2}=$ $(\alpha _{1}\alpha
_{2})^{3}=(\alpha _{1}\alpha _{3})^{2}=(\alpha _{2}\alpha _{3})^{5}=$ $%
(\alpha _{1}\alpha _{2}\alpha _{3})^{5}=$ $1\ {\rm mod}\ \Gamma (11)$, the
defining relations of $A_{5}$). $A_{5}$ is a (maximal) subgroup of $%
L_{2}(11) $. For both fixing groups we use the notation $\tilde{\Gamma}%
_{h,g}=11||11+$ \cite{FMN}. All the singular cusps of the GMF are identified
under the corresponding $\Gamma _{h,g}$ which is therefore a genus zero
fixing group.

The $q$ expansions to $O(q^{5})$ of the hauptmoduls for $\tilde{\Gamma}%
_{h,g}=11||11+$ are \cite{Q1},\cite{Q2}: 
\begin{eqnarray}
Z\left[ 
\begin{array}{c}
h \\ 
g
\end{array}
\right] (11\tau ) &=&\frac{1}{q}+0+(\frac{1}{2}\pm i\frac{\sqrt{11}}{2}%
)q+2q^{2}+(\frac{1}{2}\pm i\frac{\sqrt{11}}{2})q^{3}  \nonumber \\
&&-(1\pm i\sqrt{11})q^{4}-(\frac{1}{2}\pm i\frac{\sqrt{11}}{2})q^{5}+...
\label{11||11+ haupt}
\end{eqnarray}

Assuming that the GMF is replicable the head character expansion for $h\in
G_{11+}\equiv \mathrm{M}_{12}$ in terms of the irreducible characters $\chi
_{i}(h)$ of $\mathrm{M}_{12}$ (in ATLAS notation \cite{CCNPW}) is \cite{Q1},%
\cite{Q2} 
\begin{eqnarray}
Z\left[ 
\begin{array}{c}
h \\ 
g
\end{array}
\right] (11\tau ) &=&\frac{1}{q}+0+(\chi _{1}(h)+\chi _{4}(h))q+(\chi
_{1}(h)+\chi _{6}(h))q^{2}+  \nonumber \\
&&(\chi _{1}(h)+\chi _{4}(h)+\chi _{6}(h)+\chi _{7}(h))q^{3}+  \nonumber \\
&&(\chi _{1}(h)+2\chi _{5}(h)+\chi _{6}(h)+\chi _{7}(h)+\chi _{13}(h))q^{4}+
\nonumber \\
&&(2\chi _{1}(h)+2\chi _{4}(h)+\chi _{5}(h)+2\chi _{6}(h)+2\chi _{7}(h)+ 
\nonumber \\
&&\chi _{11}(h)+\chi _{12}(h)+\chi _{13}(h))q^{5}+...  \label{M12 char exp}
\end{eqnarray}

The upper signs in (\ref{11||11+ haupt}) according to (\ref{M12 char exp})
are for $h$ an element of $11A$ class of $\mathrm{M}_{12}$, the lower - for $%
h$ an element of $11B$ class of $\mathrm{M}_{12}$.$\qed $

Propositions \ref{order_h=5}-\ref{h order 11} lead us to the following
statement:

\begin{cor}
For $g=p+$, $o(h)=p$ and $p=5,7,11$ the irrational GMFs obey
\[
Z\left[ 
\begin{array}{c}
h \\ 
g
\end{array}
\right] =Z\left[ 
\begin{array}{c}
h^{d} \\ 
g^{a}
\end{array}
\right] \text{ iff }\left( \frac{ad}{p}\right) =1\text{,}
\]
where the fixing group $\widetilde{\Gamma }_{h,g}$ contains $\Gamma
_{0}(p^{2})$.
\end{cor}

\subsubsection{Irrational GMFs for $g=13+$ and $o(h)=13$}

The $\mathrm{L}_{3}\mathrm{(3)}$ group has four conjugate irrational classes
of order $13$. There are four subsets of powers of $h$ (\ref{h-a sets})
conjugate in $G_{13+}=\mathrm{L}_{3}\mathrm{(3)}$ determined by $\chi
_{D}^{(4)}\left( a\right) \in \{\pm 1,\pm i\}$. Since $g_{0}\stackrel{\bf M%
}{\sim }g_{0}^{s}$ ($s=1$,$2$,...,$p-1$) for any Monster element $g_{0}$, we
have the following four disjoint sets of conjugate elements in $\bf{M}$
defined by: 
\begin{eqnarray*}
\mathcal{S}_{0} &\equiv &\{(gh^{n})^{s}|\chi _{D}^{(4)}\left( n\right) =-i,%
\text{ }s\neq 0\ {\rm mod}\ 13\} \\
\mathcal{S}_{1} &\equiv &\{(gh^{n})^{s}|\chi _{D}^{(4)}\left( n\right) =1,%
\text{ }s\neq 0\ {\rm mod}\ 13\} \\
\mathcal{S}_{2} &\equiv &\{(gh^{n})^{s}|\chi _{D}^{(4)}\left( n\right) =i,%
\text{ }s\neq 0\ {\rm mod}\ 13\} \\
\mathcal{S}_{3} &\equiv &\{(gh^{n})^{s}|\chi _{D}^{(4)}\left( n\right) =-1,%
\text{ }s\neq 0\ {\rm mod}\ 13\} 
\end{eqnarray*}
where $\cup _{r=0}^{3}\mathcal{S}_{r}=\{g^{A}h^{B}|A,B\neq 0\ {\rm mod}\ 13\}
$.

The fixing group for the irreducible characters of $\mathrm{L}_{3}\mathrm{(3)%
}$ is $\langle *3\rangle $, so we expect invariance under $\Delta \equiv
\delta _{13}^{2}=\left( 
\begin{array}{cc}
3 & 0 \\ 
0 & 9
\end{array}
\right) {\rm mod}\ \Gamma (13)$, which is of order $3$ in $\Gamma (13)$,
and therefore the order of $\Gamma _{h,g}/\Gamma (13)$ is divisible by $3$.
The only possibility is $\Gamma _{h,g}/\Gamma (13)=A_{4}$, a maximal
subgroup of $L_{2}(13)$. Then $\Gamma _{h,g}$ contains an element $\alpha $,
such that $\Delta ^{3}=\alpha ^{2}=(\Delta \alpha )^{3}=1$ mod $\Gamma
(13)$. There are four solutions of the above relation, $\alpha (r)=\left( 
\begin{array}{cc}
11 & 2^{3-r} \\ 
92.2^{r} & 67
\end{array}
\right) $, $r=0,1,2,3$, corresponding to the fixing groups related to the
four irrational classes. The result can be rigorously formulated in the
following

\begin{prop}
For $g=13+$ and $h$ of order $13$, $h\in G_{13+}\equiv \mathrm{L}_{3}\mathrm{%
(3)}$ the class structure with $g$ Fricke, $h$ non-Fricke, one of the sets $%
\mathcal{S}_{r}$ ($r=0,1,2,3$) Fricke and the others non-Fricke gives rise
to an irrational GMF with a genus zero fixing group $\Gamma _{h,g}=\langle
\Gamma (13),\Delta ,\alpha (r)\rangle $, such that $\Gamma _{h,g}/\Gamma
(13)=A_{4}$.
\end{prop}

\textbf{Proof:} Let us assume that the elements of $\mathcal{S}_{1}$ are
Fricke, the others are non-Fricke. Firstly we show invariance under $\Delta $%
. According to (\ref{Delta like inv}) for some $d$, $Z\left[ 
\begin{array}{c}
h \\ 
g
\end{array}
\right] =Z\left[ 
\begin{array}{c}
h^{d} \\ 
g^{3}
\end{array}
\right] $. A $T^{-1}$ transformation gives us $Z\left[ 
\begin{array}{c}
gh \\ 
g
\end{array}
\right] =Z\left[ 
\begin{array}{c}
g^{3}h^{d} \\ 
g^{3}
\end{array}
\right] $ and hence $g^{3}h^{d}\in \mathcal{S}_{1}$ since $gh\in \mathcal{S}%
_{1}$ is Fricke. Hence $\chi _{D}^{(4)}(d)=1$ so $Z\left[ 
\begin{array}{c}
h \\ 
g
\end{array}
\right] (\tau )=Z\left[ 
\begin{array}{c}
h^{9} \\ 
g^{3}
\end{array}
\right] (\tau )\equiv Z\left[ 
\begin{array}{c}
h \\ 
g
\end{array}
\right] (\Delta \tau )$.

Now we demonstrate $\alpha (1)$ invariance. To this end, as before, we
consider a general transformation to the $(gh)^{11}$ Fricke twisted sector.
It is clear that $Z\left[ 
\begin{array}{c}
g^{a}h^{b} \\ 
g^{11}h^{11}
\end{array}
\right] $ cannot be rational for $(a,b)\neq (0,0)$ mod $13$. If we
choose $a$,$b$ such that $g^{a}h^{b}\in G_{g^{11}h^{11}}$ ($\phi
_{g^{11}h^{11}}(g^{a}h^{b})=1$), this GMF must be still irrational, since
there are no rational characters of order $13$. So that $Z\left[ 
\begin{array}{c}
g^{a}h^{b} \\ 
g^{11}h^{11}
\end{array}
\right] $ is either equal to $Z\left[ 
\begin{array}{c}
h \\ 
g
\end{array}
\right] $ or one of its algebraic conjugates. We may then further restrict
the choice of $a,b$ so that 
\[
Z\left[ 
\begin{array}{c}
h \\ 
g
\end{array}
\right] =Z\left[ 
\begin{array}{c}
g^{a}h^{b} \\ 
g^{11}h^{11}
\end{array}
\right] . 
\]
Taking a $T^{-n}$ transformation we have $Z\left[ 
\begin{array}{c}
g^{n}h \\ 
g
\end{array}
\right] =Z\left[ 
\begin{array}{c}
g^{a-2n}h^{b-2n} \\ 
g^{11}h^{11}
\end{array}
\right] $. For $\chi _{D}^{(4)}(n)=1$, $g^{n}h\in \mathcal{S}_{1}$ is Fricke
and hence $g^{a-2n}h^{b-2n}\in \mathcal{S}_{1}$ also, which leads to $%
(a,b)=(3,5)$, $(9,2)$, $(1,6)$. The solution $(9,2)$ provides the modular
transformation $\alpha (1)$. Thus $\Gamma _{h,g}=$ $\langle \Gamma
(11),\Delta ,\alpha (1)\rangle $. All the singular cusps of the GMF are
identified under $\Gamma _{h,g}$ which is therefore a genus zero fixing
group. Note that $\widetilde{\Gamma }_{h,g}$ does not contain $\Gamma
_{0}(169)$. The fixing groups for the other cases follow from (\ref{Gamma*k
f.gr.}).$\qed $

\section{Conclusions}

We have shown how Irrational Generalised Moonshine can be understood from
the analysis of the class structure for a pair of Monster elements of prime
order $p$. This follows from constraints originating from the abelian
orbifolding of the Moonshine Module and properties of centraliser
irreducible characters. We have explicitly demonstrated the genus zero
property for Irrational Generalised Moonshine Functions (GMFs) in all cases.
The methods developed in this paper and in \cite{IT} can in principle be
extended to analyse other GMFs towards proving the genus zero property in
general.

\section{Acknowledgements}

We are very grateful to S. Norton for providing us with information on
replicable series and to G. Mason for many valuable discussions. We
acknowledge funding from Enterprise Ireland under the Basic Research Grant
Scheme.

\end{document}